\documentclass[pdftex,12pt,a4paper,reqno]{amsart}\usepackage[T1]{fontenc}
\usepackage[latin1]{inputenc}
\usepackage[english]{babel}
\usepackage{lmodern}
\usepackage{amsthm}
\usepackage{amsmath}
\usepackage{amsfonts}
\usepackage{amssymb}
\usepackage[all]{xy}
\usepackage{calrsfs}
\usepackage{slashed}
\usepackage{babel}
\usepackage{array}
\newtheorem*{theo*}{Theorem}%[subsection]

\newcommand{\R}{\mathbb{R}}

\newcommand{\Sph}{\mathbb{S}}
\newcommand{\Pro}{\mathbb{P}}
\usepackage[usestackEOL]{stackengine}
\usepackage{etoolbox}
\patchcmd{\thmhead}{(#3)}{#3}{}{}
\usepackage[left=2cm,right=2cm,top=2cm,bottom=2cm]{geometry}
\usepackage{mdframed}

\title{}
\author{}

\everymath{\displaystyle\everymath{}}

\begin{document}

%\maketitle

%\setcounter{tocdepth}{2}

%\tableofcontents

%\newpage

\title{Conformal stereographic projections of sphere quotients are Majumdar-Papapetrou manifolds}
\author{Gautier Dietrich}
\thanks{The author was supported in part by the grant ANR-17-CE40-0034 of the French National Research Agency ANR (project CCEM)}
\address{Institut Montpelli\'erain Alexander Grothendieck\\ UMR 5149\\ Universit\'e de Montpellier\\ Place Eug\`ene Bataillon\\ 34090 Montpellier\\ France}
\address{Lyc\'ee polyvalent Bellevue\\ 135 route de Narbonne\\ 31031 Toulouse cedex 4\\ France} 
\email{gautier.dietrich@ac-toulouse.fr}
\date{}

\begin{abstract}
In this short note, we compute the conformal stereographic projection of a standard sphere quotient metric. The result is a Majumdar-Papapetrou metric, which might be useful.
\end{abstract}

\maketitle

\section*{Introduction}

In \cite{BN04}, in order to compute the Yamabe invariant of  $\R\Pro^3$, H. Bray and A. Neves notice that the conformal stereographic projection of $\R\Pro^3$ endowed with its standard metric is $\R^3-B_1(0)$ endowed with the Schwarzschild metric. We extend this result to general sphere quotients:

\begin{theo*}
The covering space of the conformal stereographic projection of the standard sphere $\Sph^n$ quotiented by a discrete subgroup $\Gamma$ of $\mathrm{Isom}(\Sph^n)$ is $\R^n$ endowed with a Majumdar-Papapetrou metric: $$g_{\text{MP},\Gamma}=\left(1+\sum_{\gamma\in\Gamma^*}\frac{m_\gamma}{|\cdot-\tilde{p}_\gamma|^{n-2}}\right)^\frac{4}{n-2}\delta_{\text{eucl}},$$ where $m_\gamma>0$ and $\tilde{p}_\gamma\in\R^n$ are explicitely known.
\end{theo*}

Bray and Neves then apply inverse mean curvature flow techniques on the boundary of $\R^3-B_1(0)$. Since, in our case, the corresponding boundary is non-connected, we do not expect that the rest of the proof can be generalized. We nevertheless hope that this computation, that we have not found in the literature, will turn out to be useful.

\section{Generalities}

Let us consider the standard sphere $(\Sph^n,\delta)$ embedded in $(\R^{n+1},\delta_\text{eucl})$. Let $p\in\Sph^n$. Let $\sigma_p$ be the stereographic projection, relatively to $p$, of $\Sph^n$ on the hyperplan $H_p:=(Op)^\perp$ of $\R^{n+1}$: $$\sigma_p : \Sph^n-\lbrace p \rbrace \longrightarrow H_p.$$

 Let $\Gamma$ be a discrete subgroup of $\mathrm{Isom}(\Sph^n)$ and $\delta_\Gamma$ be the metric induced  by $\delta$ on $\Sph^n/\Gamma$. The Green function $G_{p,\Gamma}$ of the conformal Laplacian $L_{\delta_\Gamma}$ at point $p$ is given by \cite{HJ99,Hab00}: $$\forall q\in\Sph^n- \Gamma p,\quad G_{p,\Gamma}(q)=\sum_{\gamma\in\Gamma}\frac{1}{|q-p_\gamma|^{n-2}},$$ where $|.|$ is the Euclidean norm in $\R^{n+1}$ and $p_\gamma:=\gamma(p)$.
 
Note that the metric induced on $H_p$ by $\sigma_p$ is $G_{p,\lbrace 1\rbrace}^\frac{4}{n-2}\delta$. Therefore, the Green function $\tilde{G}$ on $H_p$ verifies: $$\forall (q,r)\in (\Sph^n-\lbrace p \rbrace)^2,\quad \tilde{G}_{q,\lbrace 1\rbrace}\left(\sigma_p(r)\right)=\frac{1}{G_{p,\lbrace 1\rbrace}(q)G_{p,\lbrace 1\rbrace}(r)}G_{q,\lbrace 1\rbrace}(r),$$ which is to say:

\begin{equation}\label{distances}
\frac{1}{|\sigma_p(r)-\sigma_p(q)|}=\frac{|q-p|\cdot|r-p|}{|r-q|}.
\end{equation}

\section{Computation}

The computation of the metric induced by $\sigma_p$ on $\Sph^n/ \Gamma$ goes as follows:
\begin{displaymath}
\begin{aligned}
G_{p,\Gamma}^\frac{4}{n-2}\delta_\Gamma & = \left(\sum_{\gamma\in\Gamma}\frac{1}{|\cdot-p_\gamma|^{n-2}}\right)^\frac{4}{n-2}\delta_\Gamma \\
& = \left(1+\sum_{\gamma\in\Gamma^*}\left(\frac{|\cdot-p|}{|\cdot-p_\gamma|}\right)^{n-2}\right)^\frac{4}{n-2}G_{p,\lbrace 1\rbrace}^\frac{4}{n-2}\delta_\Gamma.
\end{aligned}
\end{displaymath}
Note that $G_{p,\lbrace 1\rbrace}^\frac{4}{n-2}\delta_\Gamma = \sigma_p^*\tilde{\delta}_\Gamma$, where $\tilde{\delta}_\Gamma$ is the metric induced by $\delta_\text{eucl}$ on $H_p/\Gamma$. Now, formula (\ref{distances}) gives: $$\forall q\in\Sph^n- \Gamma p,\quad \frac{|q-p|}{|q-p_\gamma|}=\frac{1}{|p-p_\gamma|}\cdot\frac{1}{|\sigma_p(q)-\sigma_p(p_\gamma)|}.$$
Consequently,
\begin{displaymath}
\begin{aligned}
G_{p,\Gamma}^\frac{4}{n-2}\delta_\Gamma & = \left(1+\sum_{\gamma\in\Gamma^*}\left(\frac{1}{|p-p_\gamma|}\cdot\frac{1}{|\sigma_p(\cdot)-\sigma_p(p_\gamma)|}\right)^{n-2}\right)^\frac{4}{n-2}\sigma_p^*\tilde{\delta}_\Gamma \\
& = \left(1+\sum_{\gamma\in\Gamma^*}\frac{m_\gamma}{d_\gamma^{n-2}}\right)^\frac{4}{n-2}\sigma_p^*\tilde{\delta}_\Gamma,
\end{aligned}
\end{displaymath}
where $m_\gamma:=\frac{1}{|p-p_\gamma|^{n-2}}$ and $d_\gamma:=|\sigma_p(\cdot)-\sigma_p(p_\gamma)|$.
We recover the mass of the manifold: $$m_{\delta_\Gamma}:=\lim_p G_{p,\Gamma}-\frac{1}{|\cdot- p|^{n-2}}=\sum_{\gamma\in\Gamma^*}m_\gamma.$$
The covering space of the conformal stereographic projection of $\Sph^n/\Gamma$ relatively to $p$ is therefore a Majumdar-Papapetrou manifold $(\R^n,g_{\text{MP},\Gamma})$, where singularities are located at the projection of the orbit of $p$ under $\Gamma$.

More generally, given a compact manifold $(M,g)$ and $\Gamma<\mathrm{Isom}(M)$ discrete, the stereographic projection of $(M/\Gamma,g_\Gamma)$ relatively to $p\in M/\Gamma$ verifies: $$G_{p,\Gamma}^\frac{4}{n-2}g_\Gamma = \left(1+\sum_{\gamma\in\Gamma^*}G_{p,\lbrace 1\rbrace}(p_\gamma)\cdot\tilde{G}_{p_\gamma,\lbrace 1\rbrace}\right)^\frac{4}{n-2}\tilde{g}_\Gamma.$$

\section*{Acknowledgments}

The author wishes to thank Marc Herzlich for introducing him to the work of Habermann and Jost, and Julien Cortier for indicating him the paper by Bray and Neves.

\bibliographystyle{alpha}

\bibliography{../Biblio}

\end{document}